\documentclass[12pt,reqno]{amsart}

\usepackage{amssymb,amsthm}
\usepackage{amsfonts,cmtiup,comment,stmaryrd}

\usepackage{mathrsfs,cmtiup}

\makeatletter
\def\blfootnote{\xdef\@thefnmark{}\@footnotetext}
\makeatother

\newtheorem{theorem}{Theorem}[section]

\newtheorem{lemma}[theorem]{Lemma}
\newtheorem{proposition}[theorem]{Proposition}
\newtheorem{corollary}[theorem]{Corollary}

\newtheorem{hyp}[theorem]{Hypothesis}
\newtheorem{conj}[theorem]{Conjecture}

\theoremstyle{definition}

\newtheorem{remark}[theorem]{Remark}
\newtheorem*{definition*}{Definition}

\newcommand{\g}{\gamma}
\newcommand{\F}{{\Bbb F}}

\newcommand{\al }{\alpha }

\newcommand{\bt}{\begin{theorem}}
\newcommand{\et}{\end{theorem}}
\newcommand{\bc}{\begin{corollary}}
\newcommand{\ec}{\end{corollary}}
\newcommand{\bpr}{\begin{proposition}}
\newcommand{\epr}{\end{proposition}}
\newcommand{\be}{\begin{equation}}
\newcommand{\ee}{\end{equation}}
\newcommand{\bp}{\begin{proof}}
\newcommand{\ep}{\end{proof}}
\newcommand{\bconj}{\begin{conj}}
\newcommand{\econj}{\end{conj}}
\newcommand{\bl}{\begin{lemma}}
\newcommand{\el}{\end{lemma}}
\newcommand{\bh}{\begin{hyp}}
\newcommand{\eh}{\end{hyp}}
\newcommand{\br}{\begin{remark}}
\newcommand{\er}{\end{remark}}

\let\leq=\leqslant
\let\geq=\geqslant
\numberwithin{equation}{section}
\newcommand{\ed}{\end{document}}
\dedicatory{Dedicated to Efim Zelmanov on the occasion of his 60th birthday}

\begin{document}
\title{Almost Engel finite and profinite groups}

\author{E. I. Khukhro}
\address{Sobolev Institute of Mathematics, Novosibirsk, 630090, Russia, and\newline University of Lincoln, U.K.}
\email{khukhro@yahoo.co.uk}

\author{P. Shumyatsky}

\address{Department of Mathematics, University of Brasilia, DF~70910-900, Brazil}
\email{pavel@unb.br}


\begin{abstract}
Let $g$ be an element of a group $G$. For a positive integer $n$, let $E_n(g)$ be the subgroup generated by all commutators $[...[[x,g],g],\dots ,g]$ over $x\in G$, where $g$ is repeated $n$ times. We prove that if $G$ is a profinite group such that for every $g\in G$ there is $n=n(g)$ such that $E_n(g)$ is finite, then $G$ has a finite
normal subgroup $N$ such that $G/N$ is locally nilpotent. The proof uses the Wilson--Zelmanov theorem saying that Engel profinite groups are locally nilpotent. In the case of a finite group $G$, we prove that if, for some $n$, $|E_n(g)|\leq m$ for all $g\in G$, then the order of the nilpotent residual $\gamma _{\infty}(G)$ is bounded in terms of $m$.
\end{abstract}
\maketitle

\section{Introduction}

A group $G$ is called an Engel group if for every $x,g\in G$ the equation $[x,g,g,\dots , g]=1$ holds, where $g$ is repeated in the commutator sufficiently many times depending on $x$ and $g$. (Throughout the paper, we use the left-normed simple commutator notation
$[a_1,a_2,a_3,\dots ,a_r]=[...[[a_1,a_2],a_3],\dots ,a_r]$.) A group is said to be locally nilpotent if every finite subset generates a nilpotent subgroup. Clearly, any locally nilpotent group is an Engel group. Wilson and Zelmanov \cite{wi-ze} proved the converse for profinite groups: any Engel profinite group is locally nilpotent. In this paper we prove related results
in terms of the subgroups
$$
E_n(g)=\langle [x,\underbrace{g,\dots ,g}_{n}]\mid x\in G\rangle.
$$

\bt\label{t-e}
Suppose that $G$ is a profinite group such that for every $g\in G$ there is a positive integer $n=n(g)$ such that $E_{n}(g)$ is finite.
Then $G$ has a finite normal subgroup $N$ such that $G/N$ is locally nilpotent.
\et

In Theorem~\ref{t-e} it also follows that there is an open locally nilpotent subgroup (just consider $C_G( N)$) --- but this fact will actually be one of the steps in the proof. The proof uses the aforementioned Wilson--Zelmanov theorem and a similar result for finite groups. Obviously, for finite groups there must be a quantitative analogue of the hypothesis that the subgroups $E_n(g)$ are finite.
For finite groups it is convenient to introduce the subgroups $E(g)=\bigcap _{n=1}^{\infty}E_n(g)$. It is also convenient to denote the nilpotent residual of a group $G$ by $\g _{\infty}(G)=\bigcap _i\g _i(G)$, where $\g _i(G)$ are terms of the lower central series ($\g _1(G)=G$, and $\g_{i+1}(G)=[\g _i(G), G]$).

\bt\label{t-emf}
Suppose that $G$ is a finite group and there is a positive integer $m$ such that $|E(g)|\leq m$ for every $g\in G$. Then the order of the nilpotent residual $\g _{\infty }(G)$ is bounded in terms of $m$.
\et

In Theorem~\ref{t-emf} it also follows that the index of the Fitting subgroup $F(G)$ is bounded in terms of $m$ (just consider $C_G( \g _{\infty }(G)$).
Theorem~\ref{t-emf} can be viewed as a generalization of Zorn's theorem \cite[Satz~III.6.3]{hup}, which says that a finite Engel group is nilpotent.

 Theorem~\ref{t-emf} also has an immediate consequence for profinite groups if there is a uniform bound for the orders of the subgroups $E_n(g)$, with the correspondingly stronger conclusion.

\bc\label{c-em}
Suppose that $G$ is a profinite group and there is a positive integer $m$ such that for every $g\in G$
 there is $n=n(g)$ such that $|E_{n}(g)|\leq m$.
Then $G$ has a finite normal subgroup $N$ of order bounded in terms of $m$ such that $G/N$ is locally nilpotent.
\ec

In a more general context, we mention our recent paper \cite{khu-shu152}, in which the subgroups $E_n(g)$ were introduced. These subgroups were used in \cite{khu-shu152} for obtaining generalizations of Baer's theorem \cite[Satz~III.6.15]{hup} saying that any Engel element of a finite group belongs to its Fitting subgroup. It was proved in \cite{khu-shu152} that if, in a soluble finite group $G$, a subgroup $E_n(g)$ has Fitting height $k$, then $g\in F_{k+1}(G)$, where $F_i(G)$ are terms of the Fitting series defined by induction:  $F_1(G)=F(G)$ is the Fitting subgroup and then $F_{i+1}(G)$ is  the inverse image of $F(G/F_i(G))$. For nonsoluble finite groups, it was proved that if the generalized Fitting height of $E_n(g)$ is $k$, then $g\in F^*_{f(k,m)}(G)$, where $F^*_i(G)$ are terms of the generalized Fitting series and $m$ is the number of prime divisors of $|g|$. A~similar result was also obtained in terms of the nonsoluble length of $E_n(g)$.
\medskip

We deal with finite groups first in \S\,\ref{s-f}, and then consider profinite groups in \S\,\ref{s-pf}.

Our notation and terminology is standard; for profinite groups, see, for example, \cite{wilson}.

We say for short that an element $g$ of a group $G$ is an \textit{Engel element} if for any $x\in G$ we have $[x,g,g,\dots , g]=1$, where $g$ is repeated in the commutator sufficiently many times depending on $x$ (such elements $g$ are often called left Engel elements).

For a group $A$ acting by automorphisms on a group $B$ we use the usual notation for commutators $[b,a]=b^{-1}b^a$ and $[B,A]=\langle [b,a]\mid b\in B,\;a\in A\rangle$, and for centralizers $C_B(A)=\{b\in B\mid b^a=b \text{ for all }a\in A\}$ and $C_A(B)=\{a\in A\mid b^a=b\text{ for all }b\in B\}$.

Throughout the paper we shall write, say, ``$(a,b,\dots )$-bounded'' to abbreviate
``bounded above in terms of $a, b,\dots
$ only''.

\section{Finite almost Engel groups}
\label{s-f}

First we list a few elementary facts that will be used without special references.

 Clearly, the subgroup $E_n(\bar g)$ constructed for the image $\bar g$ of an element $g\in G$ in a quotient $G/N$ is the image of $E_n(g)$, and if $g\in H\leq G$, then $E_n(g)$ constructed for $H$ is contained in $E_n(g)$ constructed for $G$. It follows that the condition that $|E(g)|\leq m$ for all $g\in G$ is inherited by any section of $G$.

 The following are well-known
properties of coprime actions: if $\al $ is an automorphism of a finite
group $G$ of coprime order, $(|\al |,|G|)=1$, then $ C_{G/N}(\al )=C_G(\al )N/N$ for any $\al $-invariant normal subgroup $N$, the equality $[G,\al ]=[[G,\al ],\al ]$ holds, and if $G$ is in addition abelian, then $G=[G,\al ]\times C_G(\al )$.

 In a finite group $G$ the nilpotent residual subgroup $\g _{\infty}(G)$ is of course equal to some subgroup $\g _n(G)$ for (all) sufficiently large $n$.  Clearly, $\gamma _{\infty}(G)N/N = \gamma _{\infty}(G/N)$ for any normal subgroup $N $.
 Recall that the Fitting series starts with the Fitting subgroup $F_1(G)=F(G)$, and by induction, $F_{k+1}(G)$ is the inverse image of $F(G/F_k(G))$.  The following lemma is well known and is easy to prove (see, for example, \cite[Lemma~10]{khu-maz}).

 \bl\label{l-metan}
 If $G$ is a finite group of Fitting height 2, then $\g _{\infty}(G)=\prod _q [F_q,G_{q'}]$, where $F_q$ is a Sylow $q$-subgroup of $F(G)$, and $G_{q'}$ is a Hall ${q'}$-subgroup of $G$.
 \el

\bl\label{l0}
If $P$ is a finite $p$-group, and $g\in G$ is a $p'$-element, then $[P,g]\leq E(g)$.
\el

\bp
For the abelian $p$-group $V=[P,g]/\Phi ([P,g])$ we have $V=[V,g]$ and $C_V(g)=1$. Then $V=\{[v,g]\mid v\in V\}$ and therefore also
$$
V=\{[v,\underbrace{g,\dots ,g}_n\,]\mid v\in V\}
$$
for any $n$. Hence, $E(g)\Phi ([P,g])\geq [P,g]$, whence the result.
\ep

\bl\label{l2}
Let $V$ be an elementary abelian $q$-group, and $U$ a $q'$-group of
automorphisms of $V$. If $|[V,u]|\leq m$ for every $u\in U$, then $|[V,U]|$ is $m$-bounded, and therefore $|U|$ is also $m$-bounded.
\el

\bp
First suppose that $U$ is abelian. We consider $V$ as an $\F _qU$-module. Pick $u_1\in U$ such that $[V,u_1]\ne 0$. By Maschke's theorem, $V= [V,u_1]\oplus C_V(u_1)$, and both summands are $U$-invariant, since $U$ is abelian. If $C_U([V,u_1])=1$, then $|U|$ is $m$-bounded and $[V,U]$ has $m$-bounded order being generated by $[V,u]$, $u\in U$. Otherwise pick $1\ne u_2\in C_U([V,u_1])$; then $V= [V,u_1] \oplus [V,u_2] \oplus C_V(\langle u_1,u_2\rangle )$. If $1\ne u_3\in C_U([V,u_1]\oplus [V,u_2])$, then $V= [V,u_1]\oplus [V,u_2]\oplus [V,u_3] \oplus C_V(\langle u_1,u_2,u_3\rangle )$, and so on. If $C_U([V,u_1]\oplus\dots \oplus [V,u_k])=1$ at some $m$-bounded step $k$, then again $[V,U]$ has $m$-bounded order. However, if there are too many steps, then
for the element $w=u_1u_2\cdots u_k$ we shall have $0\ne [V,u_i]= [[V,u_i],w]$, so that $[V,w] = [V,u_1]\oplus\dots \oplus [V,u_k]$ will have order greater than $m$, a contradiction.

We now consider the general case. Since every element $u\in U$ acts faithfully on $[V,u]$, the exponent of $U$ is $m$-bounded. If $P$ is a Sylow $p$-subgroup of $U$, let $M$ be a maximal normal abelian subgroup of $P$. By the above, $|[V,M]|$ is $m$-bounded. Since $M$ acts faithfully on $[V,M]$, we obtain that $|M|$ is $m$-bounded. Hence $|P|$ is $m$-bounded, since $C_P(M)= M$ and $P/M$ embeds in the automorphism group of $M$.  Since $|U|$ has only $m$-boundedly many prime divisors, it follows that $|U|$ is $m$-bounded. Since $[V,U]=\sum_{u\in U}[V,u]$, we obtain that $|[V,U]|$ is also $m$-bounded.
\ep

\bl\label{l3}
If $G$ is a finite group such that $|E(g)|\leq m$ for all $g\in G$, then $G/F(G)$ has $m$-bounded exponent.
\el

\bp
For every $g\in G$, the subgroup $E(g^k)$ is $g$-invariant for any positive integer $k$. Choose $k$ to be the maximum exponent of $\operatorname{Aut}H$ over all groups $H$ of order at most $m$. Clearly, $k$ is $m$-bounded. Then $[E(g^k), g^k]=1$. This implies that $g^k$ is an Engel element and therefore belongs to the Fitting subgroup $F(G)$ by Baer's theorem \cite[Satz~III.6.15]{hup}.
\ep

\bp[Proof of Theorem~\ref{t-emf}] Recall that $G$ is a finite group such that
$|E(g)|\leq m$ for every $g\in G$. We need to show that $|\g _{\infty }(G)|$ is $m$-bounded.

First suppose that $G$ is soluble. Since $G/F(G)$ has $m$-bounded exponent by Lemma~\ref{l3}, the Fitting height of $G$ is $m$-bounded, which follows from the Hall--Higman theorems \cite{ha-hi}. Hence we can use induction on the Fitting height, with trivial base when the group is nilpotent and $\g _{\infty }(G)=1$. When the Fitting height is at least 2, consider the second Fitting subgroup $F_2(G)$. By Lemma \ref{l-metan} we have $\g _{\infty }(F_2(G))=\prod _q [F_q,H_{q'}]$, where $F_q$ is a Sylow $q$-subgroup of $F(G)$, and $H_{q'}$ is a Hall ${q'}$-subgroup of $F_2(G)$, the product taken over prime divisors of $|F(G)|$. For a given $q$, let $\bar H_{q'}=H_{q'}/C_{H_{q'}}(F_q)$, and let $V$ be the Frattini quotient $F_q/\Phi (F_q)$. Note that $\bar H_{q'}$ acts faithfully on $V$, since the action is coprime \cite[Satz~III.3.18]{hup}.

For every $x\in \bar H_{q'}$ we have $|[V,x]|\leq m$ because $[V,x]$ is a section of $E(x)$ by Lemma~\ref{l0}. Then $|\bar H_{q'}|$ is $m$-bounded by Lemma~\ref{l2}. As a result, $|[F_q , H_{q'}]|= |[F_q ,\bar H_{q'}]|$ is $m$-bounded, since $[F_q ,\bar H_{q'}]$ is the product of $m$-boundedly many subgroups $[F_q ,\bar h]$ for $h\in H_{q'}$, for each of which $|[F_q ,h]|\leq m$, since $[F_q ,h]\leq E(h)$ by Lemma~\ref{l0}.

For the same reasons, the primes $q$ for which $[F_q , H_{q'}]\ne 1$ are less than or equal to $m$. As a result, $|\g _{\infty }(F_2(G))|$ is $m$-bounded. Induction on the Fitting height applied to $G/\g _{\infty }(F_2(G))$ completes the proof in the case of soluble $G$.

Now consider the general case. First we show that the quotient $G/R(G)$ by the soluble radical is of $m$-bounded order. Recall that the generalized Fitting subgroup is the product of the Fitting subgroup and all subnormal quasisimple subgroups (here a group is quasisimple if it is equal to its derived subgroup and its quotient by the centre is a non-abelian simple group); the generalized Fitting subgroup contains its centralizer. Let $E$ be the generalized Fitting subgroup of $G/R(G)$. It suffices to show that $E$ has $m$-bounded order. Since we are considering the quotient by the soluble radical, $E=S_1\times\dots\times S_k$ is a direct product of non-abelian finite simple groups $S_i$. Since the exponent of $G/F(G)$ is $m$-bounded by Lemma~\ref{l3}, the exponent of $E$ is also $m$-bounded. Now the classification of finite simple groups implies that every $S_i$ has $m$-bounded order, and it remains to show that the number of factors is also $m$-bounded. By Shmidt's theorem \cite[Satz~III.5.1]{hup}, every $S_i$ has a non-nilpotent soluble subgroup $R_i$, for which $\g _{\infty} (R_i)\ne 1$. We apply our  theorem, which is already proved in the case of soluble groups, to $T=R_1\times \dots \times R_k$. We obtain that $|\g _{\infty} (T)|$ is $m$-bounded, whence the number of factors is $m$-bounded.

Thus, $|G/R(G)|$ is $m$-bounded. Since $|\g _{\infty} (R(G)|$ is $m$-bounded by what was proved above in the soluble case, we can consider $G/\g _{\infty} (R(G))$ and assume that $R(G)=F(G)$ is nilpotent. Therefore $|G/F(G)|$ is $m$-bounded. We now use induction on $|G/F(G)|$.  The basis of this induction includes the trivial case $G/F(G)=1$ when $\g _{\infty}(G)=1$. But the bulk of the proof deals with the case where $G/F(G)$ is a non-abelian simple group.

 Thus, we assume that $G/F(G)$ is a non-abelian simple group of $m$-bounded order.
Let $g\in G$ be an arbitrary element. The subgroup $F(G)\langle g\rangle$ is soluble, and therefore $|\g _{\infty }(F(G)\langle g\rangle)|$ is $m$-bounded by the above. Since $\g _{\infty }(F(G)\langle g\rangle)$ is normal in $F(G)$, its normal closure $\langle \g _{\infty }(F(G)\langle g\rangle) ^G\rangle$ is a product of at most $|G/F(G)|$ conjugates, each normal in $F(G)$, and therefore has $m$-bounded order. Choose a transversal $\{t_1,\dots, t_k\}$ of $G$ modulo $F(G)$ and set
$$
K=\prod _i\langle \g _{\infty }(F(G)\langle t_i\rangle) ^G\rangle,
$$
which is a normal subgroup of $G$ of $m$-bounded order. Therefore it is sufficient to obtain an $m$-bounded estimate for $|\g _{\infty }(G/K)|$. Thus, we can assume that $K=1$. We remark that then
\be\label{e-nil}
[F(G), g, \dots , g]=1\qquad \text{for any } g\in G,
 \ee
 when $g$ is repeated sufficiently many times. Indeed, $g\in F(G)t_i$ for some $t_i$, and the subgroup $F(G)\langle t_i\rangle$ is nilpotent due to our assumption that $K=1$.

We now claim that
\be\label{e-nil2}
[F(G), G, \dots , G]=1
\ee
 if $G$ is repeated sufficiently many times. It is sufficient to prove that $[F_q, G, \dots , G]=1$ for every Sylow $q$-subgroup $F_q$ of $F(G)$. For any $q'$-element $h\in G$ we have $[Q,h]=[Q,h,h]$ and therefore $[Q,h]=1$ in view of \eqref{e-nil}. Let $H$ be the subgroup of $G$ generated by all $q'$-elements. Then $G=F_qH$ since $G/F(G)$ is simple, and $[F_q,H]=1$, so that
$$
[F_q, G, \dots , G]=[F_q, F_q, \dots , F_q]=1
$$
for a sufficiently long commutator.

We finally show that $D:=\g _{\infty }(G)$ has $m$-bounded order. First we show that $D=[D,D]$. Indeed, since $G/F(G)$ is non-abelian simple, $D$ is nonsoluble and we must have $$G=F(G)[D,D].$$
Taking repeatedly commutator with $G$ on both sides and applying \eqref{e-nil2}, we obtain $D=\g _{\infty }(G)\leq  [D,D]$, so $D=[D,D]$.

Since $F(G)\cap D$ is hypercentral in $D$ by \eqref{e-nil2} and $[D,D]=D$, it follows that
$F(G)\cap D\leq Z(D)\cap [D,D]$ by the well-known Gr\"un lemma \cite[Satz~4]{gr}. Thus, $D$ is a central covering of the simple group $D/(F(G)\cap D)\cong G/F(G)$, and therefore by Schur's theorem \cite[Hauptsatz~V.23.5]{hup} the order of $D$ is bounded in terms of the $m$-bounded order of $G/F(G)$.

We now finish the proof of Theorem~\ref{t-emf} by induction on the $m$-bounded order  $k=|G/F(G)|$ proving that $|\g _{\infty }(G)|$ is $(m,k)$-bounded. The basis of this induction is the case of  $G/F(G)$ being simple: nonabelian simple was considered above, and simple of prime order is covered by the soluble case. Now suppose that $G/F(G)$ has  a nontrivial proper normal subgroup with full inverse image $N$, so that $F(G)<N\lhd G$.  Since $F(N)=F(G)$, by induction applied to $N$, the order $|\g _{\infty }(N)|$ is bounded in terms of $m$ and $|N/F(G)|<k$. Since $N/\g _{\infty }(N)\leq F( G/\g _{\infty }(N))$,  by induction applied to $G/\g _{\infty }(N)$ the order $| \g _{\infty }(G/\g _{\infty }(N) )|$ is bounded in terms of $m$ and $|G/N|<k$. As a result, $|\g _{\infty }(G)| $ is $(m,k)$-bounded, as required.
\ep

\section{Profinite almost Engel groups}\label{s-pf}

In what follows, unless stated otherwise, a subgroup of a profinite group will always mean a closed subgroup, all homomorphisms will be continuous, and quotients will be by closed normal subgroups. Of course, any finite subgroup is automatically closed. We also say that a subgroup is generated by a subset $X$ if it is generated by $X$ as a topological group.

In this section we prove Theorem~\ref{t-e}. It is convenient to state its hypothesis as follows.

\begin{hyp}\label{hyp}
For every element $g$ of a group $G$ there is a positive integer $n=n(g)$ such that $E_{n}(g)$ is finite.
\end{hyp}

Recall that pro-(finite nilpotent) groups, that is, inverse limits of finite nilpotent groups, are called \textit{pro\-nil\-po\-tent} groups.

\bl\label{l-p-n}
A profinite group satisfying Hypothesis~\ref{hyp} is pro\-nil\-po\-tent if and only if it is locally nilpotent.
\el

\bp
Of course, any locally nilpotent profinite group is pro\-nil\-po\-tent.
Conversely, suppose that $G$ is a pro\-nil\-po\-tent group satisfying Hypothesis~\ref{hyp}. Since $E_{n(g)}(g)$ is finite, 
 for any $g\in G$ there is an open normal subgroup $N$ with nilpotent quotient $G/N$ such that $E_{n(g)}(g)\cap N=1$.  The image $\bar g$ of $g$ in the nilpotent group $G/N$ is an Engel element. Since $E_{n(g)}(g)\cap N=1$, the element $g$ is an Engel element. Thus, all elements of $G$ are Engel elements, that is, $G$ is an Engel profinite group. Then $G$ is locally nilpotent by the Wilson--Zelmanov theorem \cite[Theorem~5]{wi-ze}.
\ep

 Recall that the pro\-nil\-po\-tent residual of a profinite group $G$ is $\g _{\infty}(G)=\bigcap _i\g _i(G)$, where $\g _i(G)$ are (the closures of) the terms of the lower central series; this is the smallest normal subgroup with pro\-nil\-po\-tent quotient. The following lemma is well known and is easy to prove. Here, element orders are understood as Steinitz numbers. The same results also hold in the special case of finite groups.

 \bl\label{l-res}
 {\rm (a)} The pro\-nil\-po\-tent residual $\g _{\infty}(G)$ of a profinite group $G$ is equal to the subgroup generated by all commutators $[x,y]$, where $x,y$ are elements of coprime orders.

 {\rm (b)} For any normal subgroup $N$ of a profinite group $G$ we have $\g _{\infty}(G/N)=
 \g _{\infty}(G)N/N$.
 \el

 \bp
 Part (a) follows from the characterization of pro\-nil\-po\-tent groups as profinite groups all of whose Sylow subgroups are normal. Part (b) follows from the fact that for any elements $\bar x, \bar y$ of coprime orders in a quotient $G/N$ of a profinite group $G$ one can find pre-images $x,y\in G$ which also have coprime orders.
 \ep

Before embarking on the proof of Theorem~\ref{t-e}, we prove the following generalization of Hall's nilpotency theorem \cite{hall58}, which will be used later and which we state in a stronger form than required here, in two versions, for abstract and for profinite groups.
We denote the derived subgroup of a group $B$ by $B'$.

\bpr\label{p-hall}

{\rm (a)} Suppose that $B$ is a normal subgroup of a group $A$ such that $B$ is nilpotent of class $c$ and $\g _{d}(A/B')$ is finite of order $k$.
Then the subgroup $C=C_A(\g _{d}(A/B'))=\{a\in A\mid [\g _{d}(A), a]\leq B'\}$ has finite $k$-bounded index and is nilpotent of $(c,d)$-bounded class.

{\rm (b)} Suppose that $B$ is a normal subgroup of a profinite group $A$ such that $B$ is pro\-nil\-po\-tent and $\g _{\infty}(A/B')$ is finite.
Then the subgroup $D=C_A(\g _{\infty }(A/B'))=\{a\in A\mid [\g _{\infty }(A), a]\leq B'\}$ is open and pro\-nil\-po\-tent.
\epr

\bp
(a)
Since $A/C$ embeds into $\operatorname{Aut}\g _{d}(A/B')$, the order of $A/C$ is $k$-bounded. We claim that $C$ is nilpotent of $(c,d)$-bounded class. Indeed, using simple-commutator notation for subgroups, we have
$$
[\underbrace{C,\dots ,C}_{d+1}, C, C,\dots ]\leq [[\g _d(A), C], C, \dots ]\leq [[B,B], C, \dots ],
$$
since $[\g _d(A), C]\leq B'$ by construction.
Applying repeatedly the Three Subgroup Lemma, we obtain
\begin{align*}
 [[B,B], \underbrace{C, \dots ,C}_{2d-1}, C, \dots ]&\leq \prod _{i+j=2d-1}[[B, \underbrace{C, \dots ,C}_{i}], [B,\underbrace{C, \dots ,C}_{j}[], C, \dots ]\\
 &\leq
 [[[B,\underbrace{C, \dots ,C}_{d}], B], C, \dots ]\\
 &\leq [[[B,B],B],C,\dots ].
 \end{align*}
Thus, $\g _{d+1}(C)\leq \g _2(B)$, then $\g _{(d+1)+(2d-1)}(C)\leq \g _3(B)$, then a similar calculation gives \allowbreak $\g _{(d+1)+(2d-1)+(3d-2)}(C)\leq \g _4(B)$, and so on. An easy induction shows that $\g _f(c,d)(C)\leq \g _{c+1}(B)=1$ for $1+f(c,d)=1+d(c(c+1)/2-c(c-1)/2$, so that $C$ is nilpotent of class $ f(c,d)$.

(b) As a centralizer of a normal section, $D$ is a closed normal subgroup. Since $A/D$ embeds into $\operatorname{Aut}\g _{\infty }(A/B')$, the subgroup $D$ has finite index; thus, $D$ is an open subgroup. We now show that the image of $D=C_A(\g _{\infty }(A/B'))$ in any finite quotient $\bar A$ of $A$ is nilpotent. Let bars denote the images in $\bar A$. Then $\g _{\infty }(\bar A/\bar B')=\overline{\g _{\infty }(A/B')}$ by Lemma~\ref{l-res}(b). Therefore, $\bar D\leq C_{\bar A}(\g _{\infty }(\bar A/\bar B'))$. In a finite group, $\g _{\infty }(\bar A/\bar B')=\g _{d}(\bar A/\bar B')$ for some positive integer $d$. Hence $\bar D$ is nilpotent by part (a).
\ep

Recall that in Theorem~\ref{t-e}, $G$ is a profinite group satisfying Hypothesis~\ref{hyp}, and we need to find a finite normal subgroup such that the quotient is locally nilpotent. The first step is to prove the existence of an open locally nilpotent subgroup.

\bpr\label{p-pf1}
If $G$ is a profinite group satisfying Hypothesis~\ref{hyp}, then it has an open
normal pro\-nil\-po\-tent
subgroup.
\epr

Of course, the subgroup in question will also be locally nilpotent by Lemma~\ref{l-p-n}; the result can also be stated as the openness of the largest normal pro\-nil\-po\-tent subgroup.

\bp
For every $g\in G$ we choose an open normal subgroup $N_g$ such that $E_{n(g)}(g)\cap N=1$. Then $g$ is an Engel element in $N_g\langle g\rangle$. By Baer's theorem \cite[Satz~III.6.15]{hup}, in every finite quotient of $N_g\langle g\rangle$ the image of $g$ belongs to the Fitting subgroup. As a result, the (closure of the) subgroup $[N_g, g]$ is pro\-nil\-po\-tent.

Let $\tilde N_g$ be the normal closure of $[N_g, g]$ in $G$. Since $[N_g, g]$ is normal in $N_g$, which has finite index, $[N_g, g]$ has only finitely many conjugates, so $\tilde N_g$ is a product of finitely many normal subgroups of $N_g$, each of which is pro\-nil\-po\-tent.
Hence, so is $\tilde N_g$. Therefore all the subgroups $\tilde N_g$ are contained in the largest normal pro\-nil\-po\-tent subgroup $K$.

It is easy to see that $G/K$ is an $FC$-group (that is, every conjugacy class is finite): indeed, every $\bar g\in G/K$ is centralized by the image of $N_g$, which has finite index in $G$. A~profinite $FC$-group has finite derived subgroup \cite[Lemma~2.6]{sha}. Hence we can choose an open subgroup of $G/K$ that has trivial intersection with the finite derived subgroup of $G/K$ and therefore is abelian; let $H$ be its full inverse image in $G$. Thus, $H$ is an open subgroup such that the derived subgroup $H'$ is contained in $K$.

We now consider the metabelian quotient $M=H/K'$, which also satisfies Hypothesis~\ref{hyp}, and temporarily use symbols $E_i(g)$ for the corresponding subgroups and elements of $M$.
For every pair of positive integers $i,j$, the set
$$
E_{i,j}=\{x\in M\mid |E_{i}(x)|\leq j\}
$$
is clearly closed. By Hypothesis~\ref{hyp} we have
$$
M=\bigcup _{i,j} E_{ij}.
$$
By the Baire category theorem, one of these sets contains an open subset; that is, there is an open subgroup $U$ and a coset $aU$ such that $aU\subseteq E_{n,m}$ for some $n, m$. In other words, $|E_n(au)|\leq m$ for all $u\in U$.

It follows that $|E_{2n+1}(u)|\leq m^2$ for any $u\in U$. Indeed, consider the subgroups
$$\begin{aligned}
E_{M',n}(a)&=\langle [x,\underbrace{a,\dots ,a}_{n}\mid x\in M'\rangle\quad \text{and}\\ \quad E_{M',n}(au)&=\langle [x,\underbrace{au,\dots ,au}_{n}\mid x\in M'\rangle ,
\end{aligned}
$$
which are contained in $E_{n}(a)$ and $E_{n}(au)$, respectively, and therefore have order at most $m$. Because $M$ is metabelian, it is easy to see that both  $E_{M',n}(a)$ and $E_{M',n}(au)$ are normal subgroups of $M$. In the quotient $\bar M=M/E_{M',n}(a) E_{M',n}(au)$, both $\bar M'\langle \bar a\rangle$ and $\bar M' \langle \bar a \bar u\rangle$ are normal nilpotent subgroups of nilpotency class at most $n$. Hence their product, which contains $\bar u$, is nilpotent of class at most $2n$ by Fitting's theorem. As a result, for any $x\in M$ we have
$$
[x,\underbrace{u,\dots ,u}_{2n+1}]\in [M', \underbrace{u,\dots ,u}_{2n}]\leq E_{M',n}(a) E_{M',n}(au),
$$
so that $E_{2n+1}(u)\leq E_{M',n}(au)E_{M',n}(a)$ and $|E_{2n+1}(u)|\leq |E_{M',n}(au)|\cdot |E_{M',n}(a)|\leq m^2$.

Thus, the corresponding subgroups $E_{2n+1}(u)$ constructed for $U$ satisfy the uniform inequality $|E_{2n+1}(u)|\leq m^2$ for all $u\in U$. The same inequality holds in every finite quotient $\bar U$ of $U$, to which we can therefore apply Theorem~\ref{t-emf}. As a result, $|\g _{\infty}(\bar U)|\leq k$ for some number $k=k(m)$ depending only on $m$. Then also $|\g _{\infty}( U)|\leq k$.

Let $W$ be the full inverse image of $U$, which is an open subgroup of $G$, and $ \Gamma$ the full inverse image of $\g _{\infty}( U)$.
Now let $F=C_W(\g _{\infty}( U))=\{w\in W\mid [\Gamma ,w]\leq K'\}$. By Proposition~\ref{p-hall}(b),
this is an open normal pro\-nil\-po\-tent subgroup, which completes the proof of Proposition~\ref{p-pf1}.
\ep

\bp[Proof of Theorem~\ref{t-e}]  Recall that $G$ is a profinite group satisfying Hypothesis~\ref{hyp}, and we need to show that $\g _{\infty}(G)$ is finite. Henceforth we denote by $F(L)$ the largest normal pro\-nil\-po\-tent subgroup of a profinite group $L$. We already know that $G$ has an open normal pro\-nil\-po\-tent
 subgroup, so that $F(G)$ is also open.

Since $G/F(G)$ is finite, we can use induction on $|G/F(G)|$.
       The basis of this induction includes the trivial case $G/F(G)=1$ when $\g _{\infty}(G)=1$. But the bulk of the proof deals with the case where $G/F(G)$ is a finite simple group.

Thus, we assume that $G/F(G)$ is a finite simple group (abelian or non-abelian).
Let $p$ be a prime divisor of $|G/F(G)|$, and $g\in G\setminus F(G)$ an element of order $p^n$, where $n$ is either a positive integer or $\infty$ (so $p^n$ is a Steinitz number). For any prime $q\ne p$, the element $g$ acts by conjugation on the Sylow $q$-subgroup $Q$ of $F(G)$ as an automorphism of order dividing $p^n$. The subgroup $[Q, g]$ is a normal subgroup of $Q$ and therefore also a normal subgroup of $F(G)$. The image of $[Q, g]$ in any finite quotient is contained in the image of $E_{n(g)}(g)$ by Lemma~\ref{l0}. Since $E_{n(g)}(g)$ is finite, it follows that $[Q, g]\leq E_{n(g)}(g)$ and $[Q, g]$ is finite.

Since $[Q, g]$ is normal in $F(G)$, its normal closure $\langle [Q, g]^G\rangle $ in $G$ is a product of finitely many conjugates and is therefore also finite. Let $R$ be the product of these closures $\langle [Q, g]^G\rangle $ over all Sylow $q$-subgroups $Q$ of $F(G)$ for $q\ne p$. Since $[Q, g]\leq E_{n(g)}(g)$, there are only finitely many primes $q$ such that $[Q,g]\ne 1$ for the Sylow $q$-subgroup $Q$ of $F(G)$. Therefore $R$ is finite, and it is sufficient to prove that $\g _{\infty }(G/R)$ is finite. Thus, we can assume that $R=1$. Note that then $[Q, g^a]=1$ for any conjugate $g^a$ of $g$ and any Sylow $q$-subgroup of $F(G)$ for $q\ne p$.

Choose a transversal $\{t_1,\dots, t_k\}$ of $G$ modulo $F(G)$.
 Let $G_1=\langle g^{t_1}, \dots ,g^{t_k}\rangle$. Clearly, $G_1F(G)/F(G)$ is generated by the conjugacy class of the image of $g$. Since $G/F(G)$ is simple, we have $G_1F(G)=G$. By our assumption, the Cartesian product $T$ of all Sylow $q$-subgroups of $F(G)$ for $q\ne p$ is centralized by all elements $g^{t_i}$. Hence,  $[G_1, T]=1$. Let $P$ be the Sylow $p$-subgroup of $F(G)$ (possibly, trivial). Then also $[PG_1, T]=1$, and therefore
 $$\g _{\infty }(G)=\g _{\infty }(G_1F(G))= \g _{\infty }(PG_1).$$
  The image of $\g _{\infty }(PG_1)\cap T$ in $G/P$ is contained both in the centre and in the derived subgroup of $PG_1/P$ and therefore is isomorphic to a subgroup of the Schur multiplier of the finite group $G/F(G)$. Since the Schur multiplier of a finite group is finite \cite[Hauptsatz~V.23.5]{hup}, we obtain that $\g _{\infty }(G)\cap T=\g _{\infty }( PG_1)\cap T$ is finite. Therefore  we can assume that $T=1$, in other words, that $F(G)$ is a $p$-group.

 If $|G/F(G)|=p$, then $G$ is a pro-$p$ group, so it is pro\-nil\-po\-tent, which means that $\g _{\infty }( G)=1$ and the proof is complete. If $G/F(G)$ is a non-abelian simple group, then we choose another prime $r\ne p$ dividing $|G/F(G)|$ and repeat the same arguments as above with $r$ in place of $p$. As a result, we reduce the proof to the case $F(G)=1$, where the result is obvious.

We now finish the proof of Theorem~\ref{t-e} by induction on  $|G/F(G)|$. The basis of this induction where $G/F(G)$ is a simple group was proved above. Now suppose that $G/F(G)$ has  a nontrivial proper normal subgroup with full inverse image $N$, so that $F(G)<N\lhd G$. Since $F(N)=F(G)$, by induction applied to $N$ the group $\g _{\infty }(N)$ is finite. Since $N/\g _{\infty }(N)\leq F( G/\g _{\infty }(N))$,  by induction applied to $G/\g _{\infty }(N)$ the group $ \g _{\infty }(G/\g _{\infty }(N) )$ is also finite. As a result, $\g _{\infty }(G) $ is finite, as required.
\ep

\section*{Acknowledgements}
The first
author was supported  by the Russian Science Foundation, project no. 14-21-00065,
and the second
by the Conselho Nacional de Desenvolvimento Cient\'{\i}fico e Tecnol\'ogico (CNPq), Brazil. The first author thanks  CNPq-Brazil and the University of Brasilia for support and hospitality that he enjoyed during his visits to Brasilia.

The authors thank the referee for a number of helpful comments.

\end{document}